\documentclass{ip-journal}
\usepackage{graphicx}
\usepackage{array}
\usepackage{amsmath}
\usepackage{amssymb}
\usepackage{amsthm}
\usepackage{amsfonts}
\usepackage{mathrsfs}
\usepackage{mathtools}
\usepackage{tikz-cd}
\usepackage{subcaption}
\usepackage{enumerate}
\usepackage{esint}
\usepackage{mathabx}
\usepackage{enumerate}
\usepackage{multirow}

\newcommand*\diff{\mathop{}\!\mathrm{d}}

\newtheorem{definition}{Definition}
\newtheorem{example}{Example}
\newtheorem{remark}{Remark}
\newtheorem{proposition}{Proposition}
\newtheorem{lemma}{Lemma}
\newtheorem{corollary}{Corollary}

\begin{document}

\title[Trajectory Surfaces of Framed Curvature Flow]{Trajectory Surfaces \\ of Framed Curvature Flow}

\author{Ji{\v r}{\' i} Minar{\v c}{\' i}k \and Michal Bene{\v s}}

\date{\today}

%%%%%%%%%%%%%%%%%%%%%%%%%%%%%%%%%%%%%%%%%%%%%%%%%%%%%%%%%%%%%%%%%%%%%%%

\maketitle

\begin{abstract}
This work introduces the framed curvature flow, a generalization of both the curve shortening flow and the vortex filament equation.  Here, the magnitude of the velocity vector is still determined by the curvature, but its direction is given by an associated time-dependent moving frame. After establishing local existence and global estimates, we analyze the trajectory surfaces generated by different variations of this flow, specifically those leading to surfaces of constant mean or Gaussian curvature.
\end{abstract}

\setcounter{tocdepth}{2}
\tableofcontents
\bigskip

%%%%%%%%%%%%%%%%%%%%%%%%%%%%%%%%%%%%%%%%%%%%%%%%%%%%%%%%%%%%%%%%%%%%%%%
%%%%%%%%%%%%%%%%%%%%%%%%%%%%%%%%%%%%%%%%%%%%%%%%%%%%%%%%%%%%%%%%%%%%%%%
\section{Introduction} 
\label{sec:introduction}

Curvature driven geometric flows have been extensively studied for their favorable properties and various applications across multiple pure and applied fields. We aim to take advantage of these benefits by keeping the magnitude of the local velocity equal to curvature, but at the same time expand and generalise this family of flows by letting the velocity direction be dictated by an associated time-dependent moving frame. We refer to this new class of geometric flows as the \emph{framed curvature flow}.

Although our formulation is based on the Frenet frame, the velocity vector is well defined even in the presence of vanishing curvature, where the normal and binormal vectors are undefined. In the language of Definition 2.2 from \cite{tangent_turning_sign}, the framed curvature flow is not a Frenet frame dependent geometric flow. In this way, it is a modification of the minimal surface generating flow from \cite{old_msg}, which is defined only when the torsion and curvature are positive along the whole curve. Another advantage over \cite{old_msg} is the rich configuration space enabled by the additional degrees of freedom from the moving the velocity direction field. To demonstrate the expressivity of this approach we derive variation of the flow that trace out various surfaces of interest in the latter part of the paper.

In this work, we formulate the coupled dynamics of the moving frame and the curvature driven motion, establish local existence and uniqueness for a simplified case of this motion law, provide useful global estimates for geometric and topological quantities, classify possible singularities formed during the flow and analyse generated trajectory surfaces.

The paper is organised as follows: Section~\ref{sec:introduction} introduces the framed curvature flow and prepares the notation and lemmas required for further analysis of the flow and the trajectory surfaces it produces. The analysis is divided into Sections~\ref{sec:local_analysis} and~\ref{sec:global_analysis}. While the former deals with local behaviour including the local existence and formation of singularities, the latter (Section~\ref{sec:global_analysis}) focuses on long-term behaviour by means of length and area estimates and explores the effects of the moving frame topology. Section~\ref{sec:generated_surfaces} then showcases interesting examples of flows from the configuration space of the framed curvature flow framework. Specifically, we explore flows leading to trajectory surfaces of constant mean and Gaussian curvature.

%%%%%%%%%%%%%%%%%%%%%%%%%%%%%%%%%%%%%%%%%%%%%%%%%%%%%%%%%%%%%%%%%%%%%%%
\subsection{Motivation}  
\label{subsec:motivation}

A surprising number of natural and artificial phenomena can be modeled by a one-dimensional filament in three-dimensional Euclidean space, moving according to laws expressed as partial differential equations that depend on both the environment and the shape of the filament. Simplifying complex three-dimensional dynamical systems into a moving space curve allows for faster and more scalable numerical simulations, often revealing new insights and intuitive explanations.

Examples of natural phenomena described via motion laws of space curves include the dynamics of dislocation loops in crystalline materials \cite{Mura1987,crossslip}, the motion of scroll waves in excitable media \cite{Keener1988,Sutcliffe16}, evolution of vortex filaments in liquids via the localized induction approximation \cite{Ricca1991} or quantum vortices in superfluid media using the Gross–Pitaevskii equation \cite{zuccher_ricca_22,barenghi1997}.

Evolving space curves can be used for modeling the motion of magnetic field lines in the solar corona \cite{mhd_filaments,padilla_plasma}, the dynamics of elastic rods to model hair strands in graphics \cite{Wardetzky2008,rods_new} or defects in smectic liquid crystals \cite{randall_kamien}. Geometric flows of curves have also been applied in the context of image processing \cite{immage_processing}, quantum field theory \cite{quantum_field_theory}, origami folding \cite{origami1,origami2}, cellular automata \cite{cellular_automata}, architecture \cite{remesikova}, and medicine \cite{colonoscopy}.

A significant subset of classical geometric flows can be formulated as a gradient flow of a suitable geometric energy functional \cite{masato_notes}. Most important example is the mean curvature flow, or curve shortening flow, which minimizes area, or length. It is useful to consider geometric motion laws defined as gradient flow of various other functionals. For instance, one can use the M{\" o}bius energy or other O'Hara type energies to find optimal embedings of knots \cite{blatt2016gradient,ohara_energy,cantarella_circle}. Repulsive energies, such as the tangent point energy, can be used in numerous applications in computational geometry and computer graphics \cite{repulsive_surfaces,repulsive_curves}. Elastic energies for curves lead to elasticae curves and the Willmore energy for surfaces has been used for finding the optimal torus shape, named Clifford torus \cite{clifford_torus}, and can lead to visualy appealing solutions for the sphere eversion problem by starting the evolution from the half-way models such as the Boy's surface. Furthermore, the optimal shapes with respect to the Willmore energy with additional constraints due to Canham and Helfrich leads to shapes of biological membranes found in nature \cite{brazda_calc_var}.

Besides the applications in science and engineering, various geometric flows have proven to be remarkably useful tools in theoretical fields ranging from geometrical measure theory to differential topology, enabling the proofs of many long-standing problems. This includes, but is not limited to, the use of the inverse mean curvature flow for the proof of the Penrose inequality in \cite{huisken_penrose_ineq} or the Perelman's work on Ricci flow with surgeries \cite{perelman_1,perelman_2} leading to the proof of Poincare and Geometrization conjecture. Or more recent work on Ricci flow leading to results such as the Generalized Smale conjecture \cite{gen_smale_conjecture} and the Differentiable Sphere Theorem \cite{diff_sphere_thm}. This area is still ripe for new results, particularly in the case of higher codimension motion, which typically receives less attention. For example, open problems from \cite{ghomilist} may be within the reach, provided that further analysis of framed curvature flow is pursued.

%%%%%%%%%%%%%%%%%%%%%%%%%%%%%%%%%%%%%%%%%%%%%%%%%%%%%%%%%%%%%%%%%%%%%%%
\subsection{Framed Curvature Flow}
\label{subsec:framed_curvature_flow}

Consider a family of closed curves $\{\Gamma_t\}_{t \in [0, \underline{t})}$ evolving in the time interval $[0, \underline{t})$, where $\underline{t} > 0$ is the terminal time. For a given time $t \in [0, \underline{t})$, the curve $\Gamma_t$ is represented by a parametrization $\gamma(t, \cdot) \colon S^1 \rightarrow \mathbb{R}^3$, where $S^1 = \mathbb{R} / 2\pi \mathbb{Z}$ is the unit circle. We use the standard Frenet frame notation, where ${T}$, ${N}$, and ${B}$ denote the tangent, normal, and binormal vectors, respectively. The curvature and torsion, given by the Frenet-Serret formulae, are denoted by $\kappa$ and $\tau$, respectively. Furthermore, $g := \Vert \partial_u \gamma \Vert$ is the local rate of parametrization and $\diff s = g \diff u$ is the arclength element.

There are many ways to frame a curve \cite{bishop_frame}. The Frenet frame is in some sense canonical and easy to work with, but is ill-defined at points of vanishing curvature. We define a time-dependent moving frame that is derived from the Frenet frame using an angle functional $\theta$. The normal vector associated with this moving frame will determine the direction of the velocity vector during the framed curvature flow.

\begin{definition}[$\theta$-frame]
    For an evolving curve $\{\Gamma_t\}_{t\in[0, \underline{t})}$ and a functional
    \begin{align*}
        \theta \in \mathcal{C}^{1,2}([0, \underline{t}) \times S^1; S^1),
    \end{align*}
    we define a \emph{$\theta$-frame} of $\Gamma_t$, with \emph{$\theta$-normal} $\nu_{\theta}$ and \emph{$\theta$-binormal} $\beta_{\theta}$, using a one-parameter group of rotation $\{ \mathscr{R}_{\theta} \}$, as
    \begin{align*}
        \begin{bmatrix} \nu_{\theta} \\ \beta_{\theta} \end{bmatrix} = \mathscr{R}_{\theta} \begin{bmatrix} N \\ B \end{bmatrix}
        \mbox{, where }
        \mathscr{R}_{\theta} := \begin{bmatrix} \cos \theta & \sin \theta \\ -\sin \theta & \cos \theta \end{bmatrix} \in \mathrm{SO}(2).
    \end{align*}
    We denote the moving framed curves as $\{(\Gamma_t, \theta_t)\}_{t\in[0, \underline{t})}$, where $\theta_t := \theta(t, \cdot)$.
\end{definition}
\noindent Note that $\beta_{\theta} = T \times \nu_{\theta}$ and Frenet-Serret type formulae for the $\theta$-frame read
\begin{align*}
    \partial_s
    \begin{bmatrix} 
        T \\
        \nu_{\theta} \\ 
        \beta_{\theta}
    \end{bmatrix} &= 
    \begin{bmatrix} 
        0 & \psi_1 & -\psi_2 \\
        -\psi_1 & 0 & \psi_3 \\
        \psi_2 & -\psi_3 & 0
    \end{bmatrix}
    \begin{bmatrix} 
        T \\
        \nu_{\theta} \\ 
        \beta_{\theta} 
    \end{bmatrix},
    &&
    \begin{matrix*}[l]
        \psi_1 := \kappa\cos\theta, \\
        \psi_2 := \kappa\sin\theta, \\
        \psi_3 := \tau + \partial_s \theta.
    \end{matrix*}
\end{align*}
In the context of the trajectory surface defined in Subsection~\ref{subsec:trajectory_surfaces}, $\psi_1$ and $\psi_2$ can be interpreted as the geodesic and normal curvatures of $\Gamma_t$ immersed in $\Sigma_t$, respectively. Hereafter, we refer to $\psi_3$ as the generalised torsion.

\begin{definition}[Framed curvature flow] 
    The family of evolving framed curves $\{(\Gamma_t, \theta_t )\}_{t \in [0, \underline{t})}$ is a solution to the \emph{framed curvature flow} if the parametrization $\gamma$ and the angle functional $\theta$ satisfy the following initial-value problem
    \begin{subequations}
        \label{eq:framed-curvature-flow}
        \begin{align}
            \partial_t \gamma &= \kappa \nu_{\theta} & \partial_t \theta &= \upsilon_{\theta}
            & \mbox{in }[0, \underline{t}) \times S^1, \\
            \gamma\vert_{t=0} &= \gamma_0 & \theta\vert_{t=0} &= \theta_0
            & \mbox{in }S^1,
        \end{align}
    \end{subequations}
    where $\gamma_0$ and $\theta_0$ are the initial conditions and the $\theta$-velocity
    \begin{align*}
        \upsilon_{\theta} \in \mathcal{C}^1 ([0, \underline{t}) \times S^1; \mathbb{R})
    \end{align*}
    will be specified later in Subsections~\ref{subsec:local_existence}, \ref{subsec:cmc} and \ref{subsec:cgc}.
\end{definition}

\begin{figure}
    \centering
    \begin{subfigure}{0.5\textwidth}
        \centering
        \vspace{.4cm}
        \includegraphics[scale=0.73]{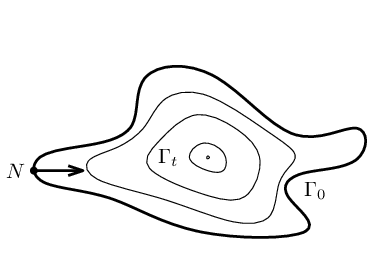}
        \caption{Curve shortening flow.}
        \label{subfig:csf}
    \end{subfigure}%
    \begin{subfigure}{0.5\textwidth}
        \centering
        \includegraphics[scale=0.73]{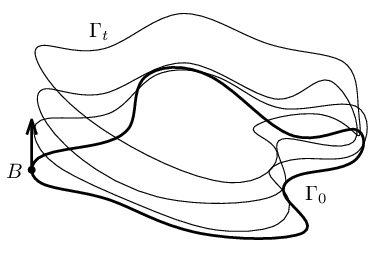}
        \caption{Vortex filament equation.}
        \label{subfig:vfe}
    \end{subfigure}
    
    \caption{Classical examples of space curve motion laws.}
    \label{fig:csf_vfe}
\end{figure}

\begin{example}\label{ex:csf_vfe}
    The framed curvature flow subsumes the following classical geometric flows, depicted in Figure~\ref{fig:csf_vfe}, as its special cases:
    \begin{enumerate}[(a)]
        \item \emph{Curve shortening flow} studied e.g. in \cite{altschuler-grayson,altschuler}: \\
        When $\theta\vert_{t = 0} = 0$ and $\upsilon_{\theta}\vert_{\theta = 0} = 0$, \eqref{eq:framed-curvature-flow} reduces to $\partial_t \gamma = \kappa N$.
        
        \item \emph{Vortex filament equation} studied e.g. in \cite{Ricca1991}: \\
        When $\theta\vert_{t = 0} = \frac{\pi}{2}$ and $\upsilon_{\theta}\vert_{\theta = \frac{\pi}{2}} = 0$, \eqref{eq:framed-curvature-flow} reduces to $\partial_t \gamma = \kappa B$.

    \end{enumerate}
\end{example}

\begin{remark}\label{rem:frenet_frame_independent}
    The framed curvature flow \eqref{eq:framed-curvature-flow} can also be viewed as a local harmonic combination of the \emph{curve shortening flow} and the \emph{vortex filament equation} from Examples~\ref{ex:csf_vfe}(a) and \ref{ex:csf_vfe}(b), resp. One can also write \eqref{eq:cmc-motion-law} as
    \begin{align*}
        \partial_t \gamma = \cos\theta \; \partial^2_s \gamma + \sin\theta \; \partial_s \gamma \times \partial^2_s \gamma.
    \end{align*}
    This formulation makes clear how the framed curvature flow \eqref{eq:framed-curvature-flow} is well-defined even when the curvature vanishes and the Frenet frame is undefined.
\end{remark}

\begin{remark}
    The set of equations \eqref{eq:framed-curvature-flow} represents a case of geometric motion with an additional quantity, namely $\theta$, whose velocity depends on the geometry, and vice versa.  This kind of coupling has been studied in e.g. \cite{Chandeliers}, where the the additional quantity represents the local radius of a bubble vortex tube. % and general relativity?
\end{remark}

%%%%%%%%%%%%%%%%%%%%%%%%%%%%%%%%%%%%%%%%%%%%%%%%%%%%%%%%%%%%%%%%%%%%%%%
%%%%%%%%%%%%%%%%%%%%%%%%%%%%%%%%%%%%%%%%%%%%%%%%%%%%%%%%%%%%%%%%%%%%%%%
\subsection{Trajectory Surfaces}
\label{subsec:trajectory_surfaces}

Similar to how a point mass moving in a homogeneous gravitational field generates a parabola, trajectory surfaces are generated by geometric flows of space curves. As the title suggests, these surfaces are one of the primary concerns of this paper. 

We argue that there are two main benefits to examining trajectory surfaces. First, the shape of the trajectory surface encodes the long-term properties of the associated motion law, and thus the knowledge of the generated surfaces may help us understand the overall behaviour of the original geometric flow. Conversely, this framework provides an alternative way to generate and study surfaces with prescribed characteristics, potentially enabling new ways to categorize and understand these surfaces and possibly help tackle various open problems. The formal meaning of trajectory surface is clarified below.

\begin{definition}[Trajectory surface, \cite{old_msg}] For a given $\theta$-velocity $\upsilon_{\theta}$, terminal time $\underline{t}$ and initial curve $\Gamma_0$, we formally define the trajectory surface $\Sigma_{\underline{t}}$ as 
\begin{align*}
\Sigma_{\underline{t}} := {\textstyle \bigcup\limits_{t \in [0,\underline{t})}} \Gamma_t,
\end{align*}
i.e. $\Sigma_{\underline{t}}$ is a surface parametrized by $\gamma(t, u)$ for $t \in [0, \underline{t})$ and $u \in S^1$.
\end{definition}

\begin{figure}
    \centering
    \includegraphics[scale=0.42]{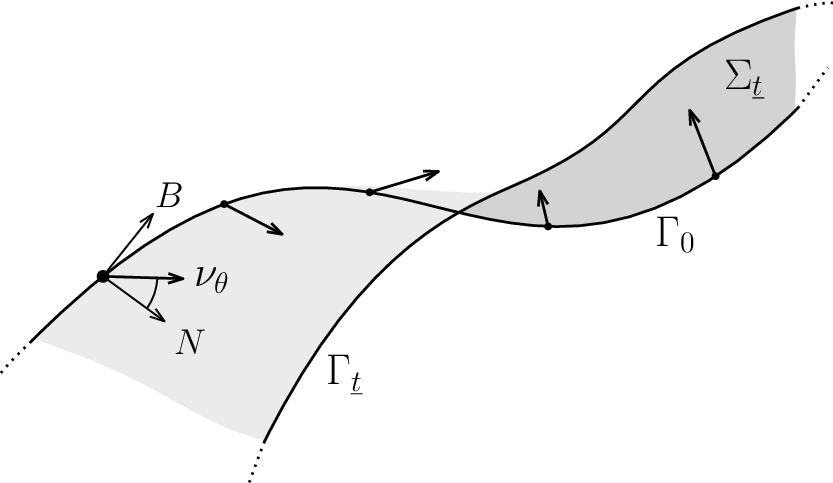}
    \caption{Trajectory surface of a framed curvature flow.}
    \label{fig:framed}
\end{figure}

\noindent Trajectory surfaces have been studied in \cite{Hussien2016} for the special case of inextensible flows, i.e. geometric flows satisfying $\partial_t g = 0$. An important example of such motion law is the vortex filament equation, mentioned in Example~\ref{ex:csf_vfe}. Surfaces generated by this motion law, referred to as Hasimoto surfaces, have been previously considered in \cite{Hussien2012}.

Closely related to the trajectory surface is the concept of worldsheet from physics. In the context of string theory, particles sweep out worldlines and strings sweep out worldsheets in Minkowski space. The equations of motion are induced from the Nambu-Goto action or the Polyakov action \cite{nambu_g_action,n_goto_action}. In our case, time is not treated as another dimension as in general relativity, but rather as another parameter.

In this paper, we are specifically interested in trajectory surfaces of constant curvature (see Section~\ref{sec:generated_surfaces}). In light of this, the following lemma states the formulas for mean and Gaussian curvature of surfaces generated by \eqref{eq:framed-curvature-flow}.

\begin{lemma}\label{lem:KH_lemma}
    Mean curvature $H$ and Gaussian curvature $K$ of the trajectory surface $\Sigma_t$ obtained from $\Gamma_t$ evolving according to \eqref{eq:framed-curvature-flow} can be expressed as
    \begin{align*}
        H = -\psi_2 + \chi \mbox{ and }  K = -\psi_3^2 - \psi_2 \chi,
    \end{align*}
    respectively. The auxiliary variable $\chi$ used in the formulae above reads
    \begin{align*}
        \chi := \frac{\upsilon_{\theta}}{\kappa} + \frac{\kappa \partial_s \psi_3 + 2 \partial_s \kappa \psi_3}{\kappa^3} \psi_1 + \frac{\partial_s^2\kappa - \kappa\partial_s \psi_3^2}{\kappa^3} \psi_2.
    \end{align*}
\begin{proof}
    The first and the second fundamental form  $\mathrm{I}$ and $\mathrm{I\!I}$ of $\Sigma_t$ read
    \begin{align*}
        \mathrm{I} &= \begin{bmatrix} 
            \mathscr{E} & \mathscr{F} \\
            \mathscr{F} & \mathscr{G} 
        \end{bmatrix} = \begin{bmatrix} 
            g_{uu} & g_{vu} \\
            g_{uv} & g_{vv}
        \end{bmatrix} = \begin{bmatrix} 
            g^2 & 0 \\
            0 & \kappa^2
        \end{bmatrix}, &
        \mathrm{I\!I} &= \begin{bmatrix} 
            \mathscr{L} & \mathscr{M} \\
            \mathscr{M} & \mathscr{N} 
        \end{bmatrix},
    \end{align*}
        where $(g_{ij})$ is the metric tensor of $\Sigma_t$, $\mathscr{L}=-g^2 \psi_2$, $\mathscr{M}= g \psi_3\kappa $ and
    \begin{align*}
        \mathscr{N} &= \kappa \upsilon_{\theta} 
        + \psi_1 \frac{\kappa\partial_s\psi_3 +2\partial_s \kappa \psi_3}{\kappa}
        + \psi_2 \frac{\partial_s^2\kappa - \kappa \psi_3^2}{\kappa}.
    \end{align*}
    Finally, the mean curvature $H$ and the Gaussian curvature $K$ are
    \begin{align*}
        K &= \frac{\det \mathrm{I\!I}}{\det \mathrm{I}} = \frac{\mathscr{L}\mathscr{N} - \mathscr{M}^2}{\mathscr{E}\mathscr{G} - \mathscr{F}^2}, &
        H &= \mathrm{tr}\left(\mathrm{I\!I}\left(\mathrm{I}^{-1}\right)\right) = \frac{\mathscr{L}}{\mathscr{E}} + \frac{\mathscr{N}}{\mathscr{G}}.
    \end{align*}    
    For more details, we refer the reader to Section 2 in \cite{old_msg}. 
\end{proof}
\end{lemma}

\begin{remark}
    The principle curvatures of the trajectory surface $\Sigma_{\underline{t}}$ generated during \eqref{eq:framed-curvature-flow} are $\kappa_{1,2} = -\psi_2 + \chi \pm \sqrt{\zeta}$, where $\zeta := \psi_2^2 - \psi_2\chi +\chi^2 +\psi_3^2$ and $\chi$ is the auxiliary variable from Lemma~\ref{lem:KH_lemma}.
\end{remark}

\noindent Further analysis of trajectory surfaces has been recently carried out in \cite{zhong2022}, which includes a description of properties of $u$-curves, e.i. curves given by $\gamma(\cdot, u)$ with a fixed parameter $u \in S^1$.

%%%%%%%%%%%%%%%%%%%%%%%%%%%%%%%%%%%%%%%%%%%%%%%%%%%%%%%%%%%%%%%%%%%%%%%
%%%%%%%%%%%%%%%%%%%%%%%%%%%%%%%%%%%%%%%%%%%%%%%%%%%%%%%%%%%%%%%%%%%%%%%
\section{Local Analysis}
\label{sec:local_analysis}

This section focuses on local properties, both in time and parameter space, of the solution to the framed curvature flow equation~\eqref{eq:framed-curvature-flow}. In particular, we state the evolution equations of the local geometric quantities in Subsection~\ref{subsec:evolution_equations}, study the effects of non-trivial tangential redistribution in Subsection~\ref{subsec:tangential} and with the help of these preliminary results we establish the local existence and uniqueness of the solution in Subsection~\ref{subsec:local_existence}. The Subsection~\ref{subsec:analysis_of_singularities} provides an overview of possible singularities formed during curvature blow-up events.

%%%%%%%%%%%%%%%%%%%%%%%%%%%%%%%%%%%%%%%%%%%%%%%%%%%%%%%%%%%%%%%%%%%%%%%
\subsection{Evolution Equations}
\label{subsec:evolution_equations}

Evolution equations for local geometric quantities, like the rate of parametrisation, curvature or torsion, during general geometric flows of space curves have been extensively studied in many pieces of literature before. See e.g. \cite{Olver08} for a general algebraic approach or \cite{sevcovic_interacting} for the treatment of geometric motion law similar to \eqref{eq:framed-curvature-flow}. Nevertheless, we state these equations and adopt them for the specific case of framed curvature flow for reader's convenience.

\begin{lemma}\label{lem:commutator}
    The arc-length commutator $[ \partial_t, \partial_s ]$ during the framed curvature flow \eqref{eq:framed-curvature-flow} with the angle functional $\theta$ is given by
    \begin{align} \label{eq:commutator}
        [ \partial_t, \partial_s] := \partial_t \partial_s - \partial_s \partial_t = \kappa^2 \cos\theta \; \partial_s.
    \end{align}
    Equivalently, $\partial_t g = -\kappa^2 \cos\theta \; g = -\kappa \psi_1 g$.
\begin{proof}
    The statement is a special case of Proposition 1 from \cite{sevcovic_interacting}.
\end{proof}
\end{lemma}

\begin{lemma}
    The Frenet frame during the framed curvature flow \eqref{eq:framed-curvature-flow} satisfies
    \begin{align*}
        \partial_t
        \begin{bmatrix} T \\ N \\ B \end{bmatrix} = 
        \begin{bmatrix} 
            0 & \xi_1 & -\xi_2 \\
            -\xi_1 & 0 & \xi_3  \\
            \xi_2 & -\xi_3  & 0
        \end{bmatrix}
        \begin{bmatrix} T \\ N \\ B \end{bmatrix},&&
        \begin{matrix*}[l]
            \xi_1 := \partial_s \psi_1 - \tau \psi_2, \\
            \xi_2 := -\partial_s \psi_2 - \tau \psi_1, \\
            \xi_3 := \kappa^{-1} (\psi_1 \partial_s \tau +\partial_s^2\psi_2 -\tau^2),
        \end{matrix*}
    \end{align*}
    while the evolution of $\theta$-normal $\nu_{\theta}$ and $\theta$-binormal $\beta_{\theta}$ can be expressed as
    \begin{align*}
        \partial_t
        \begin{bmatrix} T \\ \nu_{\theta} \\ \beta_{\theta} \end{bmatrix} = 
        \begin{bmatrix} 
            0 & \zeta_1 & -\zeta_2 \\
            -\zeta_1 & 0 & \zeta_3  \\
            \zeta_2 & -\zeta_3  & 0
        \end{bmatrix}
        \begin{bmatrix} T \\ \nu_{\theta} \\ \beta_{\theta} \end{bmatrix},&&
        \begin{matrix*}[l]
            \zeta_1 := \partial_s \kappa, \\
            \zeta_2 := -\psi_3 \kappa, \\
            \zeta_3 := \upsilon_{\theta} + \xi_3.
        \end{matrix*}
    \end{align*}
    Finally, the curvature $\kappa$ and torsion $\tau$ evolve as
    \begin{align} \label{eq:dt_curvature}
        \partial_t \kappa &= \kappa^2 \psi_1 + \kappa^{-1} \left( \partial_s^2 \psi_1 - \partial_s \tau \psi_2 - 2 \tau \partial_s \psi_2 + \tau \psi_1 \right), \\
        \partial_t \tau &= \kappa \psi_1 (\tau + \psi_3) + \partial_s \left[ \kappa^{-2} \left( \partial_s^2 \psi_2 + 2\partial_s \kappa \psi_1\psi_3 - \kappa \psi_2\psi_3 + \kappa\psi_1 \partial_s \psi_2 \right) \right]. \label{eq:dt_torsion}
    \end{align}
\begin{proof}
    Proved by substitution to Example 5.7 from \cite{Olver08}.
\end{proof}
\end{lemma}

\noindent The evolution equations for the $\theta$-frame local quantities $\psi_1$, $\psi_2$ and $\psi_3$ are more involved, but can be expressed as
\begin{align*}
    \partial_t \psi_1 &= \partial_t\kappa \cos \theta - \psi_2 \upsilon_{\theta},\\
    \partial_t \psi_2 &= \partial_t\kappa \sin \theta + \psi_1 \upsilon_{\theta},\\
    \partial_t \psi_3 &= \partial_t \tau + \partial_s\upsilon_{\theta} + \kappa \psi_1 \partial_s \theta,
\end{align*}
where $\partial_t \kappa$ and $\partial_t \tau$ shall be substituted from \eqref{eq:dt_curvature} and \eqref{eq:dt_torsion}.

%%%%%%%%%%%%%%%%%%%%%%%%%%%%%%%%%%%%%%%%%%%%%%%%%%%%%%%%%%%%%%%%%%%%%%%
\subsection{Tangential Redistribution}
\label{subsec:tangential}

To simplify previous calculations, we ignored the tangential velocity in \eqref{eq:framed-curvature-flow} by setting $\upsilon_T := \langle \partial_t \gamma, T \rangle = 0$. Apart from advection of the $\theta$-frame along the curve, this choice does not affect the geometry of the moving curve. Non-trivial tangential velocity can, however, be useful for improving numerical stability and existence analysis. We wish to do the latter in the following subsection. Hence we introduce and analyse appropriate tangential term here. Specifically, we use the tangential velocity term developed and used in \cite{HouLowengrub,KimuraTang,MikulaSevcTang} and modify it for our motion law in the following lemma.

\begin{lemma}\label{lem:tangential}
    Assume that for all $t \in [0, \underline{t})$ and all $s \in \mathbb{R} / L(\Gamma_t) \mathbb{Z}$ the tangential velocity $\upsilon_T$ satisfies the following integro-differential equation
    \begin{align}\label{eq:tangential_condition}
        \upsilon_T(t, s) = \upsilon_{T, 0}(t) + \int_{0}^{s} \kappa \psi_1 \diff\overline{s} - \frac{s}{L(\Gamma_t)} \int_{\Gamma_t} \kappa \psi_1 \diff\overline{s},
    \end{align}
    where $\upsilon_{T, 0}(t) = \upsilon_{T}(t, 0)$ is any differentiable function $\upsilon_{T, 0} \in \mathcal{C}^1([0,\underline{t}))$. Then the quantity $L(\Gamma_t)^{-1} g$ is constant during the framed curvature flow \eqref{eq:framed-curvature-flow}.
\begin{proof}
    With $\partial_t \gamma = \kappa \nu_{\theta} + \upsilon_T T$, the arc-length commutator from Lemma~\ref{lem:commutator} is
    \begin{align*}
        [\partial_t, \partial_s] = ( \kappa \psi_1 - \partial_s \upsilon_T ) \partial_s
    \end{align*}
    and we have $\partial_t g = (-\kappa\psi_1 + \partial_s \upsilon_T)g$. Note that the choice of $\upsilon_T$ does not affect the evolution of length $L(\Gamma_t)$, provided the curves $\Gamma_t$ are closed. And thus
    \begin{align*}
        \partial_t \left( \frac{g}{L(\Gamma_t)} \right) &= \frac{g}{L(\Gamma_t)^2} \left[ (-\kappa\psi_1 +\partial_s \upsilon_T) L(\Gamma_t) + \int_{\Gamma_t} \kappa \psi_1 \diff s \right].
    \end{align*}
    Substitution of $\partial_s \upsilon_T$ from \eqref{eq:tangential_condition} yields the vanishing right-hand side. 
\end{proof}
\end{lemma}

\noindent With a suitable choice of parametrization and the tangential velocity satisfying \eqref{eq:tangential_condition} we can achieve uniform parametrization throught the flow.

\begin{proposition}\label{prop:tangential}
    Assume that the initial curve $\Gamma_0$ is uniformly parametrized such that $g(0, u) = L(\Gamma_0)$ for all $u \in S^1$. Let $\{(\Gamma_t, \theta_t)\}_{t\in[0,\underline{t})}$ be a solution to the framed curvature flow with tangential velocity \eqref{eq:tangential_condition}. Then the curve $\Gamma_t$ is parametrized uniformly during the whole flow, i.e. $g(t, u) = L(\Gamma_t)$ for all $t \in [0, \underline{t})$ and $u \in S^1$.
\begin{proof}
    Straightforward application of Lemma~\ref{lem:tangential}. 
\end{proof}
\end{proposition}

\noindent Note that $\upsilon_T$ in \eqref{eq:tangential_condition} is indeed well defined on the periodic domain as one can easy verify that $\upsilon_T\vert_{s=0} \equiv \upsilon_T\vert_{s=L(\Gamma_t)}$ and $\partial_s\upsilon_T\vert_{s=0} \equiv \partial_s\upsilon_T\vert_{s=L(\Gamma_t)}$.

%%%%%%%%%%%%%%%%%%%%%%%%%%%%%%%%%%%%%%%%%%%%%%%%%%%%%%%%%%%%%%%%%%%%%%%
\subsection{Local Existence}
\label{subsec:local_existence}

\noindent This subsection establishes local existence and uniqueness for the framed curvature flow constrained by assumptions outlined in Lemma~\ref{lem:max_princ} or Lemma~\ref{lem:max_princ2}. First, it is important to note that the right-hand side of \eqref{eq:framed-curvature-flow} is well-defined even in the absence of the Frenet frame (see Remark~\ref{rem:frenet_frame_independent}).

The existence result is achieved by extending the method of abstract theory of analytic semi-flows in Banach spaces from \cite{DaPratoGrisvard,Angenent1,Angenent2,Lunardi}. In particular we formulate \eqref{eq:framed-curvature-flow} in terms of an extended four-dimensional system by treating $\theta$ as another dimension, and follow the existence proof of a similar system of equations from \cite{sevcovic_interacting}. First, let $\hat{\gamma} \colon [0, \underline{t}) \times S^1 \rightarrow \mathbb{R}^4$ denote the extended parametrization $\hat{\gamma} := [\gamma_1, \gamma_2, \gamma_3, \theta]^T$. And consider the extended system
\begin{align}\label{eq:extended_system}
    \partial_t \hat{\gamma} = \hat{\mathcal{A}} \partial_s^2\hat{\gamma} + f(\partial_s \hat{\gamma}, \hat{\gamma}),
\end{align}
where $f \in \mathcal{C}(\mathbb{R}^{4,4}; \mathbb{R}^4)$ and the principal part of the right-hand side reads
\begin{align*}
    \hat{\mathcal{A}} = \left[
        \setlength{\arraycolsep}{1.75pt}
        \def\arraystretch{1.2}
        \begin{array}{@{}c|c@{}}
        \mathcal{A} & \begin{array}{@{}c@{}} 0 \end{array} \\
        \cline{1-1}
        \multicolumn{1}{@{}c}{
            \begin{matrix} \beta^T \end{matrix}
        } & \alpha
        \end{array}
    \right],\hspace{5pt}
    \mathcal{A} = \cos \theta \; \mathbb{I} + \sin \theta \left[ T \right]_{\times} \hspace{-2pt}=\hspace{-2pt} 
    \setlength{\arraycolsep}{2pt}
    \begin{bmatrix} 
        \setlength{\arraycolsep}{-1pt}
        \cos \theta & -\sin \theta \; T_3 & \sin \theta \; T_2\\ 
        \sin \theta \; T_3 & \cos \theta & -\sin \theta \; T_1 \\ 
        -\sin \theta \; T_2 & \sin \theta \; T_1 & \cos \theta 
    \end{bmatrix},
\end{align*}
where $\mathbb{I}_{ij} = \delta_{ij}$, $(\left[ T \right]_{\times})_{ij} = \sum_k\varepsilon_{ijk} T_k$ for $i, j \in \{ 1, 2, 3 \}$, and $\alpha \in \mathbb{R}^+$, $\beta \in \mathbb{R}^3$ are fixed parameters of the framed curvature flow \eqref{eq:framed-curvature-flow} with the following $\theta$-velocity
\begin{align}\label{eq:upsilon_theta}
    \upsilon_{\theta} = \alpha \partial_s^2 \theta + \kappa \left\langle \beta, N \right\rangle +  f_4(\partial_s \hat{\gamma}, \hat{\gamma}).
\end{align}
We want the system \eqref{eq:extended_system} to be parabolic. As the spectrum of $\hat{\mathcal{A}}$ is
\begin{align*}
    \sigma(\hat{\mathcal{A}}) = \sigma(\mathcal{A}) \cup \{ \alpha \} = \{ \alpha, \cos \theta, e^{\pm i\theta} \},
\end{align*}
the eigenvalues $\alpha$ and $\cos\theta$ must be positive. In order to proceed towards the local existence result, additional constraints have to be laid down to ensure that this property is guaranteed. In the following lemmas, we provide two different ways to achieve this goal.

\begin{lemma}\label{lem:max_princ}
    Let $\beta = 0$ and $\alpha > 0$ be fixed parameters of the extended system of equations \eqref{eq:extended_system} with $f_4 \equiv 0$ and assume that $\theta_0$ satisfies $\vert \theta(0, u) \vert < \frac{\pi}{2}$ for all $u \in S^1$. Then any solution $\{(\Gamma_t, \theta_t)\}_{t\in[0, \underline{t}]}$ to \eqref{eq:framed-curvature-flow} will satisfy $\vert \theta \vert < \frac{\pi}{2}$ everywhere and the extended system \eqref{eq:extended_system} will remain parabolic.
\begin{proof}
    The statement is a consequence of the weak maximum principle for the angle functional $\theta$. Using the notation from Chapter 7.1.4 of \cite{evans_max_princ}, we have
    \begin{align}\label{eq:parabolic_operator}
        \partial_t \theta + \mathcal{L}\theta = 0,
    \end{align}
    where $\mathcal{L} := -\alpha \partial_s^2 =- \alpha g^{-2} \partial_u^2 + \alpha g^{-3} \partial_u g \partial_u$. Thanks to the trivial right-hand side of \eqref{eq:parabolic_operator}, we can use Theorem 8 from Chapter 7 of \cite{evans_max_princ} and conclude that
    \begin{align*}
        \vert \theta(t, u) \vert \leq \max_{\overline{u} \in S^1} \vert \theta(0, \overline{u}) \vert < \frac{\pi}{2}
    \end{align*}
    for all $(t, u) \in [0, \underline{t}] \times S^1$. The last inequality holds due to assumptions. 
\end{proof}
\end{lemma}

\noindent Introducing additional assumptions on the curvature allows us to extend the result from Lemma~\ref{lem:max_princ} for the case of non-trivial $\beta$ and $f_4$ from \eqref{eq:upsilon_theta}.

\begin{lemma}\label{lem:max_princ2}
    Assume $\vert\theta_0\vert < \frac{\pi}{2}$ on $S^1$ and there exist $C_1, C_2 > 0$ such that for all $u \in S^1$ and all $t \in [0, \underline{t}]$ we have $\kappa(t, u) \leq C_1$ and $f_4(t, u) \leq C_2$, where
    \begin{align*}
        \underline{t} := \frac{1}{C_1 \vert\beta\vert + C_2} \left( \frac{\pi}{2} - \max_{u \in S^1} \vert \theta_0(u) \vert \right).
    \end{align*}
    Then $\vert \theta \vert < \frac{\pi}{2}$ holds everywhere on $[0, \underline{t}] \times S^1$ and \eqref{eq:extended_system} remains parabolic. 
\begin{proof}
    The non-difusive term $\upsilon_{\theta} - \alpha \partial_s^2 \theta$ of the equation \eqref{eq:upsilon_theta} is bounded as 
    \begin{align*}
        \vert \upsilon_{\theta} - \alpha \partial_s^2 \theta \vert \leq C_3,
    \end{align*}
    where $C_3 := C_1 \vert\beta\vert + C_2$ is a positive constant. Using this value we construct
    \begin{align*}
        \theta_{-} &:= \min_{u \in S^1} \vert \theta_{0} (u) \vert - C_3, & \theta_{+} &:= \max_{u \in S^1} \vert \theta_{0} (u) \vert + C_3,
    \end{align*}
    which are subsolution and supersolution to $\theta$ (see Lemma~\ref{lem:max_princ}). Since $\vert \theta_{\pm}\vert$ stays bellow $\frac{\pi}{2}$ for all $t \in [0, \underline{t}]$, we concur that $\vert \theta \vert$ is bounded by $\frac{\pi}{2}$ as well.
\end{proof}
\end{lemma}

\noindent To prepare for the existence proof, further notation needs to be introduced. For any $\varepsilon \in (0, 1)$ and any $k \in \{0, \frac{1}{2}, 1\}$ we define the following family of Banach spaces of H{\" o}lder continuous functions
\begin{align*}
    \mathcal{E}_{k} := h^{2k + \varepsilon}(S^1) \times h^{2k + \varepsilon}(S^1) \times h^{2k + \varepsilon}(S^1) \times h^{2k + \varepsilon}(S^1),
\end{align*}
where $h^{2k + \varepsilon}(S^1)$ is a little H{\" o}lder space (see Section 4.1 in \cite{sevcovic_interacting}). With the aid of the previous lemmas and the appropriate tangential velocity term described in Subsection~\ref{subsec:tangential}, we can now state the local existence result.

\begin{proposition}\label{prop:local_exist}
    Consider \eqref{eq:framed-curvature-flow} with additional tangential velocity satisfying \eqref{eq:tangential_condition} from in Subsection~\ref{subsec:tangential} and assume that
    \begin{enumerate}[\;\;\;(a)]
        \item the initial extended parametrization $\hat{\gamma}\vert_{t=0}$ belongs to $\mathcal{E}_1$,
        \item the initial parametrization $\gamma_0$ satisfies $\Vert \partial_u \gamma_0 \Vert = L(\Gamma_0)$ on $S^1$,
        \item $f$ is $\mathcal{C}^2$ smooth and globally Lipschitz continuous,
        \item the assumptions of Lemma~\ref{lem:max_princ} or Lemma~\ref{lem:max_princ2} are satisfied.
    \end{enumerate}
    Then there exists $\underline{t} > 0$ and a unique family of framed curves $\{(\Gamma_t, \theta_t)\}_{t\in[0, \underline{t})}$ satisfying \eqref{eq:extended_system} with tangential velocity \eqref{eq:tangential_condition} such that $\hat{\gamma} \in \mathcal{C}([0, \underline{t}); \mathcal{E}_1) \cap \mathcal{C}^1([0, \underline{t}); \mathcal{E}_0)$.
\begin{proof}
    We extend the proof of Theorem 4.1 from \cite{sevcovic_interacting} to the parametrization with the framing angle $\hat{\gamma}$. We rewrite the extended system \eqref{eq:extended_system} as an abstract parabolic equation:
    \begin{align}\label{eq:abstract_parabolic_eq}
        \partial_t \hat{\gamma} + \mathscr{F}(\hat{\gamma}) = 0
    \end{align}
    for $\hat{\gamma} \in \mathcal{E}_{1}$, where $\mathscr{F}$ is operator mapping from $\mathcal{E}_{1}$ to $\mathcal{E}_{0}$. Using Lemma 2.5 from \cite{Angenent1} as in the proof of Proposition 4 from \cite{sevcovic_interacting}, the Frechet derivative $\mathscr{F}'$ of the operator $\mathscr{F}$ from \eqref{eq:abstract_parabolic_eq} belongs to the maximum regularity class $\mathcal{M}(\mathcal{E}_1, \mathcal{E}_0)$. The solution $\hat{\gamma}$ exists in 
    \begin{align*}
        \mathcal{C}([0, \overline{t}]; \mathcal{E}_1) \cap \mathcal{C}^1([0, \overline{t}]; \mathcal{E}_0)
    \end{align*}
    for any $\overline{t} \in (0, \underline{t})$ due to Theorem 2.7 from \cite{Angenent1}. 
\end{proof}
\end{proposition}

\noindent For more details, we refer the reader to \cite{sevcovic_interacting} or to the original literature \cite{DaPratoGrisvard,Angenent1,Angenent2,Lunardi} of the abstract theory of analytic semi-flows in Banach spaces. 

\begin{proposition}\label{prop:blow_up}
    Let the assumptions of Proposition~\ref{prop:local_exist} hold and suppose that the maximal time of existence $\underline{t}$ is finite, then
    \begin{align*}
        \limsup_{t \to \underline{t}} \max_{u \in S^1} \kappa(t, u) = +\infty \;\; \vee \;\;  \limsup_{t \to \underline{t}} \max_{u \in S^1} \vert \partial_s^2 \theta(t, u) \vert = +\infty.
    \end{align*}
\begin{proof}
    For contradiction, assume that the maximal time of the existence $\underline{t}$ is finite and that both $\kappa$ and $\vert \partial_s^2 \theta \vert$ are bounded. Since the assumptions of Proposition~\ref{prop:local_exist} are satisfied, the solution $\hat{\gamma}$ belongs to $\mathcal{C}([0, \overline{t}]; \mathcal{E}_1) \cap \mathcal{C}^1([0, \overline{t}]; \mathcal{E}_0)$ for any $\overline{t} \in (0, \underline{t})$. Moreover, $\partial_s^2 \hat{\gamma}$ is bounded by the assumptions because
    \begin{align*}
        \Vert \partial_s^2 \hat{\gamma} \Vert^2 = \Vert \partial_s^2 \gamma \Vert^2 + \vert \partial_s^2 \theta \vert^2 = \kappa^2 + \vert \partial_s^2 \theta \vert^2.
    \end{align*}
    Thus, by the maximum regularity, the extended solution $\hat{\gamma}$ belongs to the space $\mathcal{C}([0, \underline{t}]; \mathcal{E}_1) \cap \mathcal{C}^1([0, \underline{t}]; \mathcal{E}_0)$ and can be continued beyond $[0, \underline{t})$, which contradicts the maximal time assumption. More details can be retrieved from the last part of Theorem 4.1 from \cite{sevcovic_interacting}. 
\end{proof}
\end{proposition}

\noindent The behaviour of $\Gamma_t$ during the curvature blow-up event described in Proposition~\ref{prop:blow_up} is detailed in the next subsection.

%%%%%%%%%%%%%%%%%%%%%%%%%%%%%%%%%%%%%%%%%%%%%%%%%%%%%%%%%%%%%%%%%%%%%%%
\subsection{Formation of Singularities}
\label{subsec:analysis_of_singularities}

The expressive power of the framed curvature flow framework allows for the occurrence of unusual singularity types during curvature blow-up events. In general, understanding these singularities has been crucial for analyzing the behavior of various geometric flows. For our work, it will be particularly important for the global analysis in Section~\ref{sec:global_analysis}.

Singularities of geometric flows have been studied in e.g. \cite{Khan15,Corr16,litzinger,andrews_singularities,ishitawa_blow_up} or \cite{altschuler-grayson}, where the motion of planar curve has been extended beyond curvature singularities via a higher dimensional flow of an associated space curve. The existence of flows past various singularities has also been addressed by other means such as by using the concept of viscosity solutions for the level-set formulation of curvature driven flows in \cite{level_set_method,pozar_viscosity}, topological surgeries \cite{perelman_2} or by analysis of self-similar shrinkers in \cite{vega_corners}.

In \cite{altschuler}, Altschuler showed that the blow-up limits of space curves under the curve shortening flow are planar. The situation for framed curvature flow is more complicated. The following definition clarifies the meaning of different types of singularities which may occur during \eqref{eq:framed-curvature-flow}.

\begin{definition}[Singularity typologies]\label{def:flat_singularity}
The event at which the curvature $\kappa$ approaches infinity at time $\underline{t}$ and $L(\Gamma_t) \rightarrow 0$ as $t$ approaches $\underline{t}$ (the curve $\Gamma_t$ shrinks to a point) during the framed curvature flow \eqref{eq:framed-curvature-flow} is called:
    \begin{enumerate}
        \item \emph{Flat singularity} if and only if $\theta_t \rightrightarrows 0$ as $t$ approaches $\underline{t}$. \\($\theta$-frame uniformly approaches the Frenet frame).
        \item \emph{Cone singularity} if and only if $\theta_t \rightrightarrows \Theta \in (-\frac{\pi}{2}, \frac{\pi}{2}) \setminus \{0\}$ as $t$ approaches $\underline{t}$.
        \item \emph{Pinch singularity} if and only if $\theta_t \rightrightarrows \Theta \in \{\pm\frac{\pi}{2}\}$ as $t$ approaches $\underline{t}$.
    \end{enumerate}
    Furthermore, when the trajectory surface $\Sigma_{\underline{t}}$ is unbounded, the case 3. is called \emph{infinite pinch singularity}.
\end{definition}

\noindent The flat singularity occurs in the classical example of curve shortening flow. For the case of simple planar curves, this singularity is guarantied by the Gage-Hamilton-Grayson  theorem \cite{GageEx,Grayson}. Embedded space curves under the curve shortening flow do not necessarily shrink to a point, however the Gage-Hamilton-Grayson theorem can be extended in the case of simple spherical curves \cite{MinKimBenSIAM}. Other singularity types from Definition~\ref{def:flat_singularity} are illustrated in the next set of analytical examples with simple evolution of circle.

\begin{figure}
    \centering
    \includegraphics[scale=0.85]{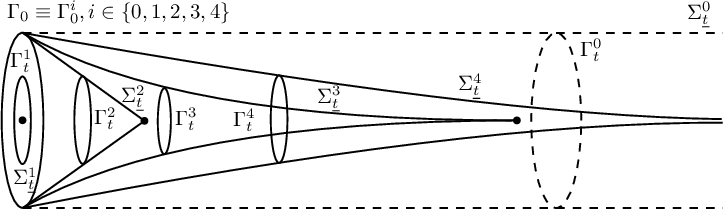}
    \caption{
        Depiction of singularity typologies from Definition~\ref{def:flat_singularity} for circular initial condition $\Gamma_0$ and uniform theta angle. Upper indices $1$ to $4$ denote the flat, cone, pinch and infinite pinch scenario, respectively. The dashed line represents the trajectory of the vortex filament equation for reference.
    }
    \label{fig:circular_singularities}
\end{figure}

\begin{example}[Cone singularity]\label{ex:cone_singularity}
    Let $\Gamma_0$ be a circle with radius $\rho_0 > 0$ and consider $\upsilon_{\theta} \equiv 0$ with $\theta_0 \equiv \phi \in (0, \frac{\pi}{2})$. This setup leads to the famous solution for shrinking circle with radius $\rho(t) = (\rho_0^2 - 2t \cos\phi)^{\frac{1}{2}}$ which vanishes at the terminal time $\underline{t} = (2\cos\phi)^{-1} \rho_0^2$. However, due to the non-trivial binormal velocity term $\langle \partial_t \gamma, B \rangle = \kappa \sin \phi \neq 0$, the singularity occurs at a point shifted in the binormal direction from the center of the initial circle by a distance
    \begin{align*}
        \underline{z} = \left(\rho_0 - \left(\rho_0^2 - 2 \underline{t} \cos\phi\right)^{\frac{1}{2}}\right) \tan\phi = \rho_0 \tan\phi.
    \end{align*}
    This leads to a conical trajectory surface with a cone singularity.
\end{example}

\begin{example}[Pinch singularity]\label{ex:pinch_singularity}
    Again, consider a circle $\Gamma_0$ with radius $\rho_0 > 0$ and $\theta_0 \equiv \phi \in (0, \frac{\pi}{2})$. To illustrate the pinch singularity, let
    \begin{align*}
        \partial_t \theta = \frac{\tan\theta - 2 \kappa \sqrt{\underline{t} - t}}{2(\underline{t} - t)},
    \end{align*}
    where $\underline{t} = (\sin\phi)^{-2}\rho_0$. This $\theta$-velocity is constructed so that $\rho = \sin\theta \sqrt{\underline{t} - t}$ and $\partial_t z = \kappa \sin\theta = \sqrt{\underline{t} - t}$. Even though the time derivative of the shift distance $z$ diverges as $t$ approaches $\underline{t}$, its definite integral remains finite:
    \begin{align*}
        \lim_{t \to \underline{t}} z(t) = \lim_{t \to \underline{t}} \int_0^t \frac{1}{\sqrt{\underline{t} - \overline{t}}} \diff \overline{t} = \lim_{t \to \underline{t}} \left[ -2 \sqrt{\underline{t} - \smash{\overline{t}}} \right]_{\overline{t}=0}^{\overline{t}=t} = 2\sqrt{\underline{t}}.
    \end{align*}
    The pinch singularity thus develops at a point located at a $2\sqrt{\underline{t}}$ distance form the center of the original circle.
\end{example}

\begin{example}[Infinite pinch singularity]\label{ex:infinite_pinch_singularity}
    With the same setup of initial conditions as in the previous examples, we now consider $\partial_t \theta = 2\rho_0^{-2}$. This velocity leads to $\theta(t) = 2t\rho_0^{-2}$ and $\rho(t) = \rho_0 (1 - \sin(2\rho_0^{-1}t))^{\frac{1}{2}}$. Thus the circle shrinks to a point at the time $\underline{t} = \tfrac{1}{4}\pi \rho_0^2$ and
    \begin{align}\label{eq:dt_z_inf_pinch}
        \partial_t z = \rho^{-1} \sin(2\rho_0^{-1}t)) = \rho_0^{-1} (1 - \sin(2\rho_0^{-1}t)))^{-\frac{1}{2}} \sin(2\rho_0^{-1}t)).
    \end{align}
    Unlike in the Example~\ref{ex:pinch_singularity}, the integral of \eqref{eq:dt_z_inf_pinch} diverges and thus the infinite pinch singularity is formed.
\end{example}

\noindent All singularity typologies from Examples~\ref{ex:cone_singularity}, \ref{ex:pinch_singularity} and \ref{ex:infinite_pinch_singularity} are depicted in Figure~\ref{fig:circular_singularities}.

\begin{remark}
    Let us recall the definition of type-I and type-II singularity, studied in e.g. \cite{Khan15,Corr16,litzinger}. This classification of curvature blow-up events is based on the comparison of the curvature growth near $\underline{t}$ with the function $(\underline{t} - t)^\frac{1}{2}$. Formally, a blow-up singularity is classified as type-I if
    \begin{align}\label{eq:singularity_type_I}
        \lim_{t \to \underline{t}} M_t (\underline{t} - t) := \lim_{t \to \underline{t}} [ \max_{u \in S^1} \kappa^2(t, u) ] (\underline{t} - t)
    \end{align}    
    is bounded, and as type-II otherwise. In terms of this classical notation, the singularities from Example~\ref{ex:cone_singularity} and \ref{ex:pinch_singularity} are type-I:
    \begin{itemize}
        \item[$\bullet$]{
            In Example~\ref{ex:cone_singularity}, we have $M_t = \rho^{-2}(t) = (\rho_0^2 - 2 t \cos \phi)^{-1} = 2 \cos \phi (\underline{t} - t)^{-1}$, where the terminal time is $\underline{t} = (2\cos\phi)^{-1} \rho_0^2$. Thus
            \begin{align*}
                \lim_{t \to \underline{t}} M_t (\underline{t} - t) = 2\cos\phi < +\infty.
            \end{align*}
        }
        \item[$\bullet$]{
            Similarly in Example~\ref{ex:pinch_singularity}, the radius reads $\rho = (\underline{t} - t)^{\frac{1}{2}} \sin \theta$ and therefore
            \begin{align*}
                \lim_{t \to \underline{t}} M_t (\underline{t} - t) = \lim_{t \to \underline{t}} \max_{u \in S^1} (\sin \theta_t(u))^{-2} = 1 < +\infty.
            \end{align*}
        }
    \end{itemize}
    Whereas the infinite pinch singularity from Example~\ref{ex:infinite_pinch_singularity} is type-II:
    \begin{itemize}
        \item[$\bullet$]{
            Since in Example~\ref{ex:infinite_pinch_singularity}, the radius is $\rho(t) = \rho_0 (1 - \sin(2\rho_0^{-1}t))^{-\frac{1}{2}}$, the term $M_t$ behaves as $(\underline{t} - t)^{2}$ near $\underline{t}$ and the limit \eqref{eq:singularity_type_I} diverges.
        }
    \end{itemize}
    Further analysis of these connections is left for a future work.
\end{remark}

\noindent The study of singularity formation is an extensive field of research, and this subsection offers only a brief exploration of potential blow-up scenarios within the context of the recently introduced framed curvature flow. Future work can involve for instance the analysis of singularity profiles leading to the self-shrinking Abresch–Langer curves \cite{abresch_langer}.

%%%%%%%%%%%%%%%%%%%%%%%%%%%%%%%%%%%%%%%%%%%%%%%%%%%%%%%%%%%%%%%%%%%%%%%
%%%%%%%%%%%%%%%%%%%%%%%%%%%%%%%%%%%%%%%%%%%%%%%%%%%%%%%%%%%%%%%%%%%%%%%
\section{Global Analysis}
\label{sec:global_analysis}

This section examines the global aspects of solutions to the framed curvature flow, with an emphasis on the properties of global geometric quantities and their long-term behavior. Subsection~\ref{subsec:global_estimates} provides several global estimates for the length and generated surface area, Subsection~\ref{subsec:projection_area} explores the evolution of the largest projected algebraic area, and Subsection~\ref{subsec:frame_topology} presents selected key facts related to the topology of the $\theta$-framing.

%%%%%%%%%%%%%%%%%%%%%%%%%%%%%%%%%%%%%%%%%%%%%%%%%%%%%%%%%%%%%%%%%%%%%%%
\subsection{Global Estimates}
\label{subsec:global_estimates}

We aim to derive useful estimates for global geometric quantities such as length and generated area. To this end, we first establish evolution equations for these quantities and then state the assumptions on which the subsequent bounds in this subsection are based.

\begin{lemma}\label{lem:dt-length-area}
    The evolution of the length $L(\Gamma_t)$ of the curve $\Gamma_t$ and the total area $A(\Sigma_t)$ of the trajectory surface $\Sigma_t$ during \eqref{eq:framed-curvature-flow} is
    \begin{align}\label{eq:dt-length-area}
        \frac{\diff}{\diff t} L ( \Gamma_t ) &= - \int_{\Gamma_t} \kappa \psi_1 \diff s, &
        \frac{\diff}{\diff t} A \left( \Sigma_t \right) &= \int_{\Gamma_{t}} \kappa \diff s.
    \end{align}
\begin{proof}
    The first part of \eqref{eq:dt-length-area} is due to \eqref{eq:commutator}, the second one follows from
    \begin{align*}
        \frac{\diff}{\diff t} \int_{\Sigma_t} \diff  A &= \frac{\diff}{\diff t} \int_{0}^{t} \int_{\Gamma_{\bar{t}}} \kappa \diff s \diff \bar{t},
    \end{align*}
    where $\diff A$ is obtained from $\mathscr{E} = g^2$, $\mathscr{F} = 0$ and $\mathscr{G} = \beta^2 + \gamma^2 = \kappa^2$ as
    \begin{align*}
        \diff  A = \sqrt{\mathscr{E}\mathscr{G} - \mathscr{F}^2} \diff u \wedge \diff t = g \kappa \diff u \wedge \diff t.
    \end{align*} 
    Particular details of the computations are in the proof of Lemma~\ref{lem:KH_lemma}.
\end{proof}
\end{lemma}

\noindent Without any assumptions on the initial curve, we can bound the generated area from below by a linear function of time.

\begin{corollary}\label{cor:area_fary}
    The Fenchel Theorem implies $A\left( \Sigma_t \right) \geq 2\pi t$ for all $t \in [0, \underline{t})$. Furthermore, when the curve is knotted on $[0, \underline{t})$, we get $A\left( \Sigma_t \right) \geq 4\pi t$ by the Milnor-F{\' a}ry Theorem \cite{Fary}.
\end{corollary}

\begin{lemma}\label{lem:dt_total_torsion}
    The evolution of the total torsion $\tau$ and the total generalized torsion $\psi_3$ of $\Gamma_t$ during the framed curvature flow \eqref{eq:framed-curvature-flow} is
    \begin{align}\label{eq:dt-total-torsion}
        \frac{\diff}{\diff t} \int_{\Gamma_t} \tau \diff s &=  \frac{\diff}{\diff t} \int_{\Gamma_t} \psi_3 \diff s = \int_{\Gamma_t} \psi_1 \psi_3 \kappa + \psi_2 \partial_s \kappa \diff s.
    \end{align}
\begin{proof}
    Since $\Gamma_t$ is closed, both integrals are equal, i.e.
    \begin{align*}
        \int_{\Gamma_t} \psi_3 \diff s = \int_{\Gamma_t} \tau \diff s + \int_{\Gamma_t} \partial_s \theta \diff s = \int_{\Gamma_t} \tau \diff s.
    \end{align*}
    The right-hand side of \eqref{eq:dt-total-torsion} is obtained from \eqref{eq:commutator} and \eqref{eq:dt_torsion}. 
\end{proof}
\end{lemma}

\noindent The estimates below are based a subset of the following \emph{assumptions}:
\begin{enumerate}[I.]
    \item There exists a fixed $\varepsilon > 0$ such that $\vert\theta\vert \leq \frac{\pi}{2} - \varepsilon$. In this case we define a constant $K_{\mathrm{I.}} := \cos\left( \tfrac{\pi}{2} - \varepsilon \right) > 0$ which will bound $\cos \theta$ from below.
    
    \item The curvature $\kappa$ is uniformly bounded from bellow by a constant $K_{\mathrm{II.}} > 0$, i.e. $\kappa(u, t) \geq K_{\mathrm{II.}}$ for all $t \in [0, \underline{t})$ and $u \in S^1$.
\end{enumerate}

\noindent Note that assumption \emph{I.} is also needed for the existence proof in Subsection~\ref{subsec:local_existence} and follows from the assumptions given in Lemma~\ref{lem:max_princ}, or alternatively Lemma~\ref{lem:max_princ2}. Assumption \emph{II.}, on the other hand, can only be enforced up-to a time $t$ away from the singularity $\underline{t}$, where the curvature blows up.

\begin{proposition}\label{prop:length_bound1}
    Let $\Gamma_t$ be a solution to \eqref{eq:framed-curvature-flow}. If assumption \emph{I.} holds, then
    \begin{align}\label{eq:length_bound1}
        L(\Gamma_t) \leq \left( L^2(\Gamma_0) - 8\pi^2 K_{\mathrm{I.}} t \right)^{\frac{1}{2}}
    \end{align}
    and thus the terminal time can be bounded from above as
    \begin{align}\label{eq:terminal_bound1}
        \underline{t} \leq \left( 8 \pi^2 K_{\mathrm{I.}} \right)^{-1} L^2(\Gamma_0).
    \end{align}
    Moreover, assuming the curve is knotted on $[0, \underline{t})$, the $8\pi^2$ term in both \eqref{eq:length_bound1} and \eqref{eq:terminal_bound1} estimates can be replaced by $32\pi^2$ as in Corollary~\ref{cor:area_fary}.
\begin{proof}
    From assumption \emph{I.} and the first part of Lemma~\ref{lem:dt-length-area} we have
    \begin{align*}
        \frac{\diff}{\diff t} L(\Gamma_t) = - \int_{\Gamma_t} \kappa^2 \cos \theta \; \diff s \leq - K_{\mathrm{I.}} \int_{\Gamma_t} \kappa^2 \diff s.
    \end{align*}
    Using the Fenchel Theorem and Cauchy-Schwarz inequality, we obtain
    \begin{align*}
        \frac{\diff}{\diff t} L(\Gamma_t) \leq - \frac{K_{\mathrm{I.}}}{L(\Gamma_t)} \left( \int_{\Gamma_t} \kappa \diff s \right)^2 \leq - \frac{4 \pi^2 K_{\mathrm{I.}}}{L(\Gamma_t)}.
    \end{align*}
    The result follows from the ODE comparison theorem.
\end{proof}
\end{proposition}

\begin{proposition}
    Let $\Gamma_t$ be a solution to \eqref{eq:framed-curvature-flow} and let $\Sigma_t$ denote its associated trajectory surface. If assumptions \emph{I.} and \emph{II.} are satisfied, then we get
    \begin{align*}
        A (\Sigma_t ) \leq \frac{ K_{\mathrm{II.}}}{ 12\pi^2 K_{\mathrm{I.}} } \left( L^3(\Gamma_0) - \left( L^2(\Gamma_0) - 8 \pi^2 K_{\mathrm{I.}} t \right)^{\tfrac{3}{2}} \right).
    \end{align*}
    Furthermore, as $\underline{t}$ is bounded by \eqref{eq:terminal_bound1}, we get a global bound
    \begin{align*}
        A (\Sigma_{\underline{t}} ) \leq \frac{ K_{\mathrm{II.}} L^3(\Gamma_0) }{ 12\pi^2 K_{\mathrm{I.}} }.
    \end{align*}
\begin{proof}
    Assuming \emph{I.} and \emph{II.} and using Lemma~\ref{lem:dt-length-area} and Proposition~\ref{prop:length_bound1} yields
    \begin{align*}
        \frac{\diff}{\diff t} A (\Sigma_t ) \leq K_{\mathrm{II.}} L(\Gamma_t) \leq K_{\mathrm{II.}} \left( L^2(\Gamma_0) - 8 \pi^2 K_{\mathrm{I.}} t \right)^{\frac{1}{2}}.
    \end{align*}
    Integrating the inequality yields the result.
\end{proof}
\end{proposition}

%%%%%%%%%%%%%%%%%%%%%%%%%%%%%%%%%%%%%%%%%%%%%%%%%%%%%%%%%%%%%%%%%%%%%%%
\subsection{Projected Area}
\label{subsec:projection_area}

In this subsection, we consider the following quantity
\begin{align}\label{eq:projected_area_integral}
    A_p(\Gamma_t) := \frac{1}{2} \int_{\Gamma_t} \gamma \times \partial_s \gamma \diff s = \frac{1}{2} \int_{S^1} \gamma \times \partial_u \gamma \diff u
\end{align}
and use it to extend our area estimates for $\Sigma_{\underline{t}}$. The geometric interpretation of this quantity is revealed in following lemma.

\begin{lemma}
    For a given curve $\Gamma_t$, the quantity  $A_p(\Gamma_t)$ defined in \eqref{eq:projected_area_integral} is the largest algebraic area enclosed by any orthogonal projection of $\Gamma_t$.
\begin{proof}
    For any unit normal vector $\nu \in S^2$, let $\pi(\nu)$ be the projection operator onto $\{ \nu \}^{\perp}$, i.e. $\pi(\nu) = \mathbb{I} - \nu \cdot \nu^T$, and let $\pi(\nu) \Gamma_t$ denote the projected planar curve paramerized by $\pi(\nu) \gamma$. Then
    \begin{align*}
        A(\pi(\nu) \Gamma_t) :=& \left\Vert \int_{S^1} \pi(\nu)\gamma \times \partial_u (\pi(\nu) \gamma) \diff u \right\Vert \\
        =& \left\Vert A_p(\Gamma_t)  - \int_{S^1} \left\langle \gamma, \nu \right\rangle \nu \times \partial_u \gamma \diff u - \int_{S^1} \left\langle \partial_u \gamma, \nu \right\rangle \gamma \times \nu \diff u \right\Vert \\
        =& \left\Vert A_p(\Gamma_t) + \nu \times \int_{S^1} \left\langle \partial_u \gamma, \nu \right\rangle \gamma - \left\langle \gamma, \nu \right\rangle \partial_u \gamma \diff u \right\Vert,
    \end{align*}
    where $A(\Gamma)$ denotes the algebraic area enclosed by a planar curve $\Gamma$. Double application of the triple vector product formula yields
    \begin{align*}
        A(\pi(\nu) \Gamma_t)
        =& \left\Vert A_p(\Gamma_t) + \nu \times (\nu \times A_p(\Gamma_t))  \right\Vert \\
        =& \left\Vert A_p(\Gamma_t) + \left\langle A_p(\Gamma_t), \nu \right\rangle \nu -  \Vert \nu \Vert^2 A_p(\Gamma_t) \right\Vert \\
        =& \left\Vert \left\langle A_p(\Gamma_t), \nu \right\rangle \nu \right\Vert = \left\vert \left\langle A_p(\Gamma_t), \nu \right\rangle \right\vert.
    \end{align*}
    Thus, due to the Cauchy-Schwarz inequality, we have:
    \begin{enumerate}[\;\;\;1.]
        \item $A(\pi(\nu) \Gamma_t) \leq A_p(\Gamma_t)$ for all $\nu$ in $S^2$,
        \item $A(\pi(A_p(\Gamma_t)^{-1} A_p(\Gamma_t)) \Gamma_t) = A_p(\Gamma_t)$ when $A_p(\Gamma_t) > 0$.
    \end{enumerate}
    The conjunction of 1. and 2. proves the statement. 
\end{proof}
\end{lemma}

\begin{lemma}\label{lem:dt_ap}
    The time derivative of $A_p(\Gamma_t)$ during \eqref{eq:framed-curvature-flow} is
    \begin{align}\label{eq:dt_ap}
        \frac{\diff}{\diff t} A_p(\Gamma_t) = - \int_{\Gamma_t} \kappa \beta_{\theta} \diff s.
    \end{align}
\begin{proof}
    The derivation of this integral quantity is simplified because
    \begin{align*}
        \int_{\Gamma_t} \gamma \times \partial_s \gamma \diff s = \int_{S^1} \gamma \times \partial_u \gamma \diff u,
    \end{align*}
    where $\diff s = g \diff u$ and, formally, $\partial_u = g \partial_s$. Thus we have
    \begin{align*}
        \frac{\diff}{\diff t} A_p(\Gamma_t) = \frac{1}{2}\int_{S^1} \partial_t \gamma \times \partial_u \gamma + \gamma \times \partial_u \partial_t \gamma \diff u = \frac{1}{2}\int_{\Gamma_t} \kappa \nu_{\theta} \times T + \gamma \times \partial_s (\kappa \nu_{\theta}) \diff s,
    \end{align*}
    where $ \nu_{\theta} \times T = -\beta_{\theta}$ and both parts of the integral yield the same value as
    \begin{align*}
        \int_{\Gamma_t} \gamma \times \partial_s (\kappa \nu_{\theta}) \diff s = \int_{\Gamma_t} \partial_s( \gamma \times \kappa \nu_{\theta} ) - T \times \kappa \nu_{\theta} \diff s = -\int_{\Gamma_t} \kappa \beta_{\theta} \diff s.
    \end{align*}
    Adding these integrals leads to \eqref{eq:dt_ap}.
\end{proof}
\end{lemma}

\begin{proposition}\label{prop:disc_projected_area}
    Let $\{(\Gamma_t, \theta_t)\}_{t \in [0, \underline{t})}$ develop a flat singularity at time $\underline{t}$, then
    \begin{align*}
        A(\Sigma_{\underline{t}}) \geq \max_{t \in [0, \underline{t})} \Vert A_p(\Gamma_t) \Vert.
    \end{align*}
\begin{proof}
    Let $x \in \mathbb{R}^3$ be the point to which the curve $\Gamma_t$ shrinks towards as the time $t$ approaches $\underline{t}$. Since $\Sigma_{\underline{t}} \cup \{ x \}$ spans the initial curve $\Gamma_0$, its area must be at least that of minimal spanning surface, which has locally larger area that the projection. Formally, let $\diff A'$ denote the area element of $\pi (\nu) \Gamma_t $, then
    \begin{align*}
        \mathscr{E}' &= \Vert \partial_u \pi(\nu) \gamma \Vert^2 \leq \vvvert \pi(\nu) \vvvert^2 \Vert \partial_u \gamma \Vert^2 = g^2,\\ 
        \mathscr{G}' &= \Vert \partial_t \pi(\nu) \gamma \Vert^2 \leq \vvvert \pi(\nu) \vvvert^2 \Vert \partial_t \gamma \Vert^2 = \kappa^2.
    \end{align*}
    With $\vvvert \cdot \vvvert$ being the spectral norm, $\vvvert \pi(\nu) \vvvert = \max(\sigma(\pi(\nu))) = 1$ and thus
    \begin{align*}
        \diff A' = \sqrt{\mathscr{E}' \mathscr{G}' - \mathscr{F}'^2} \diff u \wedge \diff t \leq \sqrt{\mathscr{E}' \mathscr{G}'} \diff u \wedge \diff t \leq g \kappa \diff u \wedge \diff t = \diff A.
    \end{align*}
    Note that the algebraic area of the projection is even smaller as the overlapping parts can annihilate. 
\end{proof}
\end{proposition}

\begin{proposition}
    Let $\Gamma_t$ be a solution to \eqref{eq:framed-curvature-flow} and assume \emph{I.} and \emph{II.}, then
    \begin{align*}
        \frac{\diff}{\diff t} \Vert A_p (\Gamma_t) \Vert \leq 2 K_{\mathrm{I.}} K_{\mathrm{II.}} L(\Gamma_t).
    \end{align*}
\begin{proof}
    Applying the assumptions and Cauchy-Schwarz inequality yields
    \begin{align*}
        \frac{\diff}{\diff t} \Vert A_p (\Gamma_t) \Vert  = \left\langle \frac{A_p(\Gamma_t)}{\Vert A_p (\Gamma_t) \Vert}, \frac{\diff}{\diff t} A_p(\Gamma_t)  \right\rangle \leq  \left\Vert \frac{\diff}{\diff t} A_p(\Gamma_t) \right\Vert \leq 2 K_{\mathrm{I.}} K_{\mathrm{II.}} \left\Vert \int_{\Gamma_t} \beta_{\theta} \diff s \right\Vert.
    \end{align*}
    The statement then follows from the fact that $\mathscr{R}_{\theta}$ is unitary. 
\end{proof}
\end{proposition}

\noindent Since the area of any surface enclosed by the curve $\Gamma_t$ is smaller than $\Vert A_p(\Gamma_t) \Vert$ (see proof of Proposition~\ref{prop:disc_projected_area}), the above proposition provides an upper bound on the growth of minimal spanning surface area.

\begin{remark}
    Note that for $\theta \equiv \frac{\pi}{2}$, both the length $L(\Gamma_t)$ and the projected area $A_p(\Gamma_t)$ remain constant during \eqref{eq:framed-curvature-flow}, as shown in \cite{arms-hama}. On the other hand, for $\theta \equiv 0$ the motion is an $L^2$-gradient flow for the length functional (see \cite{masato_notes}).
\end{remark}

%%%%%%%%%%%%%%%%%%%%%%%%%%%%%%%%%%%%%%%%%%%%%%%%%%%%%%%%%%%%%%%%%%%%%%%
\subsection{Frame Topology}
\label{subsec:frame_topology}

Let us consider the topology of the moving $\theta$-frame and its possible ramifications on the long-term behaviour of the framed curvature flow. We do so by analysing the time evolution of two geometric quantities, named writhe and twist, which are closely connected to the topology of the moving frame. Writhe of an embedded curve is an averaged sum of all signed crossings over the space of all orthogonal projections, but can also be written using Gauss formula as
\begin{align*}
    \mathrm{W}_{\mathrm{r}}(\Gamma_t) = \frac{1}{4\pi} \iint_{\Gamma_t \times \Gamma_t} \frac{\left\langle \gamma(s, t) - \gamma(s', t), T(s, t) \times T(s', t) \right\rangle}{\Vert \gamma(s, t) - \gamma(s', t) \Vert^3} \diff s \wedge \diff s'.
\end{align*}

\noindent The second important geometric quantity describing the frame topology is the total twist of the $\theta$-frame, which reads
\begin{align*}
    \mathrm{T}_{\mathrm{w}}^{\theta}(\Gamma_t) &= \frac{1}{2\pi} \int_{\Gamma_t} \langle \nu_{\theta} \times \partial_s \nu_{\theta}, T \rangle \diff s \\
    &= \frac{1}{2\pi} \int_{\Gamma_t} -\psi_1\langle \nu_{\theta} \times T, T \rangle +\psi_3\langle \nu_{\theta} \times \beta_{\theta}, T \rangle \diff s \\
    &= \frac{1}{2\pi} \int_{\Gamma_t} \psi_3 \diff s = \mathrm{T}_{\mathrm{w}}^{F} (\Gamma_t) + \frac{1}{2\pi} \mathrm{deg} (\theta_t), 
\end{align*}
where $\mathrm{deg}(\theta_t)$ is the topological degree of $\theta_t \colon S^1 \rightarrow S^1$ and $\mathrm{T}_{\mathrm{w}}^{F} (\Gamma_t)$ is the normalized total twist (i.e. total twist assiciated with the Frenet frame):
\begin{align*}
    \mathrm{T}_{\mathrm{w}}^F(\Gamma_t) = \frac{1}{2\pi} \int_{\Gamma_t} \tau \diff s.
\end{align*}

\noindent The writhe and twist are connected via the Călugăreanu–White–Fuller Theorem \cite{calugareanu,white} which states that
\begin{align}\label{eq:C-W-F}
    \mathrm{S}_{\mathrm{Lk}}^{\bullet}(\Gamma_t) &= \mathrm{W_r}(\Gamma_t) + \mathrm{T}_{\mathrm{w}}^{\bullet}(\Gamma_t),
\end{align}
where $\bullet$ represents either $F$ or $\theta$ and $\mathrm{S}_{\mathrm{Lk}}^{\theta}(\Gamma_t)$, $\mathrm{S}_{\mathrm{Lk}}^{F}(\Gamma_t)$ are the self-linking numbers for the Frenet frame and the $\theta$-frame, respectively. With the help of this theorem we can describe the evolution of writhe for embedded curves.
\begin{proposition}\label{prop:d_dt_writhe}
    Let $\{(\Gamma_t, \theta_t)\}_{t\in[0,\underline{t})}$ be a solution to \eqref{eq:framed-curvature-flow}. Consider $t \in [0, \underline{t})$ such that $\Gamma_t \hookrightarrow \mathbb{R}^3$ (i.e. $\Gamma_t$ is embedded) and $\kappa(u, t) > 0$ for all $u \in S^1$. Then
    \begin{align}\label{eq:dt_writhe}
        \frac{\mathrm{d}}{\mathrm{d} t} \mathrm{W}_{\mathrm{r}}(\Gamma_t) = - \frac{1}{2\pi} \int_{\Gamma_t} \psi_1 \psi_3 \kappa + \psi_2 \partial_s \kappa \diff s.
    \end{align}
\begin{proof}
    The assumptions imply that the time derivative of $\mathrm{S}_{\mathrm{Lk}}^F(\Gamma_t)$ exists and is equal to 0. Thus, we may differentiate \eqref{eq:C-W-F} to obtain
    \begin{align*}
        \frac{\mathrm{d}}{\mathrm{d} t} \mathrm{W}_{\mathrm{r}} (\Gamma_t) =  \frac{\mathrm{d}}{\mathrm{d} t} \left[ \mathrm{S}_{\mathrm{Lk}}^F(\Gamma_t) - \mathrm{T}_{\mathrm{w}}^F(\Gamma_t) \right] = - \frac{1}{2\pi} \frac{\mathrm{d}}{\mathrm{d} t} \int_{\Gamma_t} \tau \diff s.
    \end{align*}
    The formula \eqref{eq:dt_writhe} then follows from Lemma~\ref{lem:dt_total_torsion}. 
\end{proof}
\end{proposition}

\begin{remark}
    Note that since $\upsilon_{\theta}$ is continuous, the degree of $\theta$ cannot change during the flow and the difference $\mathrm{S}_{\mathrm{Lk}}^{\theta}(\Gamma_t) - \mathrm{S}_{\mathrm{Lk}}^{F}(\Gamma_t)$ is a constant integer.
\end{remark}

\noindent The following proposition provides a necessary topological condition needed to close the trajectory surface in sense of ending the flow in a flat singularity described in Subsection~\ref{subsec:analysis_of_singularities}.

\begin{proposition}
    Assume that  $\{(\Gamma_t, \theta_t)\}_{t \in [0, \underline{t})}$ develops a flat singularity and $(\Gamma_0, \theta_0)$ is not a Seifert framing. Then there must exist $t \in [0, \underline{t})$ such that either $\kappa(u, t) = 0$ at some point $u \in S^1$ or  $\Gamma_t$ is not embedded.
\begin{proof}
    The Seifert framing must have zero self-linking number \cite{renzo-seifert}.
\end{proof}
\end{proposition}

%%%%%%%%%%%%%%%%%%%%%%%%%%%%%%%%%%%%%%%%%%%%%%%%%%%%%%%%%%%%%%%%%%%%%%%
%%%%%%%%%%%%%%%%%%%%%%%%%%%%%%%%%%%%%%%%%%%%%%%%%%%%%%%%%%%%%%%%%%%%%%%
\section{Generated Surfaces}
\label{sec:generated_surfaces}

By adjusting the definition of the $\theta$-velocity, the framed curvature flow can be formulated such that its associated trajectory surface has various interesting properties. These specific formulations are explored in this section. We consider trajectory surfaces of constant mean curvature in Subsection~\ref{subsec:cmc} and then constant Gaussian curvature in Subsection~\ref{subsec:cgc}. Other special surfaces, such as surfaces of constant ratio of principle curvatures, proposed in \cite{new_pottman}, fall outside the scope of this manuscript and may be the subject of future work.

%%%%%%%%%%%%%%%%%%%%%%%%%%%%%%%%%%%%%%%%%%%%%%%%%%%%%%%%%%%%%%%%%%%%%%%
\subsection{Constant Mean Curvature}
\label{subsec:cmc}

In this subsection, we consider the use of framed curvature flow as a means of solving the Bj{\" o}rling problem for minimal surfaces and its generalisation for non-minimal surfaces of constant mean curvature (see \cite{bjorling}).

\begin{proposition}
    For a fixed constant $H \in \mathbb{R}$, consider the framed curvature flow \eqref{eq:framed-curvature-flow} with the $\theta$-velocity given by
    \begin{align}\label{eq:cmc-motion-law}
        \upsilon_{\theta} = \kappa H -(\kappa\partial_s\psi_3 + 2 \partial_s \kappa \psi_3) \kappa^{-2} \psi_1 + (\kappa^3 + \kappa \psi^2_3 - \partial_s^2\kappa) \kappa^{-2} \psi_2.
    \end{align}
    The trajectory surface $\Sigma_{\underline{t}}$ generated by this flow has a constant mean curvature equal to the prescribed value $H$.
\begin{proof}
    Substitution of \eqref{eq:cmc-motion-law} to Lemma~\ref{lem:KH_lemma}. 
\end{proof}
\end{proposition}

\noindent The following results are all related to the Flux theorem introduced in \cite{kusner_thesis,kusner_balanced}.

\begin{proposition}[Flux Theorem]
    Let $\{(\Gamma_t, \theta_t)\}_{t \in [0, \underline{t})}$ be a solution to the framed curvature flow with $\theta$-velocity defined in \eqref{eq:cmc-motion-law}. Then for any $a \in \mathbb{R}^3$
    \begin{align}\label{eq:flux_theorem}
        H \int_{\partial \Sigma_{t}} \langle \gamma \times T, a \rangle \diff s + \int_{\partial \Sigma_{t}} \langle \nu_{\theta}, a \rangle \diff s = 0,
    \end{align}
    where $\partial \Sigma_{t} = \Gamma_0 \cup \Gamma_t$ is the boundary of the associated trajectory surface $\Sigma_t$.
\begin{proof}
    Multiple proofs can be found in e.g. \cite{LopezBook}, where the unit conormal vector from Theorem 5.1.1. is equivalent to the $\theta$-normal $\nu_{\theta}$.
\end{proof}
\end{proposition}

\noindent Combining the Flux theorem with the evolution equations for the projected area $A_p$ leads to a simple formula for the derivative of total $\theta$-normal.

\begin{corollary}
        If $\{(\Gamma_t, \theta_t)\}_{t \in [0, \underline{t})}$ is a solution to \eqref{eq:framed-curvature-flow} with $\upsilon_{\theta}$ from \eqref{eq:cmc-motion-law}, then
    \begin{align*}
        \frac{\diff}{\diff t} \int_{\Gamma_t} \nu_{\theta} \diff s = 2 H \int_{\Gamma_t} \kappa \beta_{\theta} \diff s.
    \end{align*}
\begin{proof}
    The boundary $\partial \Sigma_t$ consists of $\Gamma_t$ and $\Gamma_0$, but the former is static. Differentiation of the Flux theorem \eqref{eq:flux_theorem} thus leads to
    \begin{align}\label{eq:dif_flux}
        H \frac{\diff}{\diff t} \int_{\Gamma_t} \gamma \times T\diff s + \frac{\diff}{\diff t} \int_{\Gamma_t} \nu_{\theta} \diff s = 0.
    \end{align}
    The result is then obtained by rearranging \eqref{eq:dif_flux} and using Lemma~\ref{lem:dt_ap}.
\end{proof}
\end{corollary}

\noindent The Flux theorem enables the following necessary conditions for flat singularity formation during the flow that generates surfaces of constant mean curvature.

\begin{corollary}
    When $\{(\Gamma_t, \theta_t)\}_{t \in [0, \underline{t})}$ develops a flat singularity at the terminal time $\underline{t}$, the surface $\Sigma_{\underline{t}}$ only has one non-trivial boundary $\Gamma_t$, and thus
    \begin{align*}
        A_p(\Gamma_0) = - \frac{1}{H} \int_{\Gamma_0} \nu_{\theta} \diff s
    \end{align*}
\end{corollary}

\begin{corollary}
    Assume that $\{(\Gamma_t, \theta_t)\}_{t \in [0, \underline{t})}$ solves the framed curvature flow with $\upsilon_{\theta}$ from \eqref{eq:cmc-motion-law} and develops a flat singularity. Then
    \begin{align*}
        L(\Gamma_0) \geq 2H \Vert A_p(\Gamma_0) \Vert.
    \end{align*}
    In particular, if $\Gamma_0$ is simple planar curve enclosing area $A$, then
    \begin{align*}
        L(\Gamma_0) \geq 2 H A.
    \end{align*}
\begin{proof}
    For any unit vector $a \in S^2$, we have $ 2H \vert \langle A_p(\Gamma_0), a \rangle \vert \leq L(\Gamma_0)$ from Corollary 5.1.7 of \cite{LopezBook}. 
\end{proof}
\end{corollary}

\noindent Important subclass of surfaces with constant mean curvature are the minimal surfaces \cite{meeks_minimal_book,golden_ages} characterised by $H=0$. In nature, minimal surfaces emerge in the context of soap films \cite{MoffattRicca10,Goldstein14}, cell membranes \cite{Larssons14}, crystallographic structure of zeolites \cite{Andersson1985,Andersson1988,Scriven1976} and as the apparent horizon of a black hole \cite{huisken_penrose_ineq}.

For the case of minimal surfaces, many of the previous results derived from the Flux theorem significantly simplify. Additionally, when the flow develops a flat singularity, the associated trajectory surface effectively represents a valid solution to the Plateau problem with single boundary curve $\Gamma_0$.

\begin{corollary}[Corollary 5.1.5 from \cite{LopezBook}]
    For minimal surfaces with $H = 0$, we have even stricter conditions, namely for all $a \in \mathbb{R}^3$
    \begin{align*}
        \int_{\partial \Sigma_t} \langle \nu_{\theta}, a \rangle \diff s &= 0, & \int_{\partial \Sigma_t}  \langle \nu_{\theta}, \gamma \times a \rangle \diff s &= 0.
    \end{align*}
\end{corollary}

\noindent We end this subsection with analysis of specific examples of solutions to \eqref{eq:cmc-motion-law} with simple initial configurations. First example illustrates the configuration that leads to the simplest minimal surface, which is a subset of the flat plane.

\begin{example}\label{ex:convex_planar}
    Let the initial curve $\Gamma_0$ be a closed convex planar curve and $\theta_0 \equiv 0$. Then the framed curvature flow with $\theta$-velocity from \eqref{eq:cmc-motion-law} and $H = 0$ is equivalent to the curve shortening flow and generates a flat surface $\Sigma_{\underline{t}}$ equivalent to the convex hull of $\Gamma_0$ in a finite time $\underline{t}$ when the $\Gamma_{t}$ shrinks down to a round point (see \cite{GageEx,Grayson}).
\end{example}

\noindent More analytical examples can be obtained by considering configurations with helical and cylindrical symmetries.

\begin{example}[Helical symmetry]\label{ex:helical_symmetry}
    For a constant $\theta_0$ consider evolving helix
    \begin{align*}
        \gamma_0(u) &:= \begin{bmatrix}  
        \varrho_0 \cos u \\
        \varrho_0 \sin u \\
        wu
        \end{bmatrix}, &
        \gamma(t, u) &:= \begin{bmatrix}  
        \varrho(t) \cos(u + \upsilon(t)) \\
        \varrho(t) \sin(u + \upsilon(t)) \\
        wu + \omega(t)
        \end{bmatrix},
    \end{align*}
    where $\rho_0$ and $w$ are positive constants and $\rho, \upsilon, \omega$ are functions of time $t$. Since $\kappa = \varrho g^{-2}$ and $\tau = w g^{-2}$, the problem \eqref{eq:cmc-motion-law} reduces to the following system
    \begin{align}\label{eq:helical_example}
        \frac{\diff}{\diff t} \begin{bmatrix} \theta \\ \varrho \\ \omega \\ \upsilon \end{bmatrix} = 
        \frac{1}{g^3}
        \begin{bmatrix}
            g \sin\theta + \rho g H \\
            - g \rho \cos\theta \\
            \rho^2 \sin\theta \\
            - w \sin\theta
        \end{bmatrix}, &&
        \left.
        \begin{bmatrix} \theta \\ \varrho \\ \omega \\ \upsilon \end{bmatrix}
        \right\vert_{t = 0}=
        \begin{bmatrix} \theta_0 \\ \varrho_0 \\ 0 \\ 0 \end{bmatrix},
    \end{align}
    where $g^2 = \varrho^2 + w^2$. Similar helicoidal surfaces of constant mean curvature were studied in e.g. \cite{helicoidal_stability}.
\end{example}

\noindent Considering cylindrically symmetrical configurations leads to the family of Delaneu surfaces, first classified in \cite{delaunay_thm}.

\begin{example}[Cylindrical symmetry]\label{ex:cylindrical_symmetry}
    Setting $w = 0$ reduces \eqref{eq:helical_example} to
    \begin{align*}
        \frac{\diff}{\diff t} \begin{bmatrix} \theta \\ \varrho \\ \omega \end{bmatrix} = 
        \frac{1}{\varrho^2}
        \begin{bmatrix}
            \sin\theta + \varrho H \\
            - \varrho \cos\theta \\
            \varrho \sin\theta
        \end{bmatrix}, &&
        \left.
        \begin{bmatrix} \theta \\ \varrho \\ \omega \end{bmatrix}
        \right\vert_{t = 0}=
        \begin{bmatrix} \theta_0 \\ \varrho_0 \\ 0 \end{bmatrix}.
    \end{align*}
\end{example}

%%%%%%%%%%%%%%%%%%%%%%%%%%%%%%%%%%%%%%%%%%%%%%%%%%%%%%%%%%%%%%%%%%%%%%%
\subsection{Constant Gaussian Curvature}
\label{subsec:cgc}

Unlike the mean curvature flow or the minimal surface generating flow \cite{old_msg}, the framed curvature flow can be used for generating developable surfaces or surfaces of any prescribed Gaussian curvature.

\begin{proposition}
    For a fixed constant $K \in \mathbb{R}$, consider a framed curvature flow \eqref{eq:framed-curvature-flow} with the $\theta$-velocity given by
    \begin{subequations}
        \label{eq:cgc-motion-law}
        \begin{align}
            \upsilon_{\theta} =& -\kappa \psi_2^{-1} K -(\kappa\partial_s\psi_3 + 2 \partial_s \kappa \psi_3) \kappa^{-2} \psi_1 \\
            &- \kappa \psi_2^{-1} \psi_3^2 - (\partial_s^2\kappa - \kappa \psi^2_3) \kappa^{-2} \psi_2.
        \end{align}
    \end{subequations}
    The trajectory surface $\Sigma_{\underline{t}}$ generated by this flow has a constant Gaussian curvature equal to the prescribed value $K$.
\begin{proof}
    Substitution of \eqref{eq:cgc-motion-law} to Lemma~\ref{lem:KH_lemma}. 
\end{proof}
\end{proposition}

\noindent Important examples of surfaces with constant Gaussian curvature are developable surfaces. For this specific case, the Gauss-Bonnet formula significantly simplifies and can be used to uncover an unexpected integral of motion.

\begin{proposition}\label{prop:cgc_developable}
    Let $\{(\Gamma_t, \theta_t)\}_{t \in [0, \underline{t})}$ be a solution to the framed curvature flow with $\theta$-velocity defined in \eqref{eq:cgc-motion-law} with $K$ set to $0$. Then the integral of $\psi_1$ over the curve $\Gamma_t$ at any time $t \in[0, \underline{t})$ is preserved.
\begin{proof}
    The Gauss-Bonnet theorem states that
    \begin{align}\label{eq:gauss_bonnet}
        \int_{\Sigma_t} K \diff A + \int_{\partial \Sigma_t} \kappa_g \diff s = 2 \pi \chi(\Sigma_t),
    \end{align}
    where $\diff A = \kappa g \diff u \wedge \diff t$ (see proof of Proposition~\ref{prop:disc_projected_area}), $\kappa_g$ is the geodesic curvature at the boundary $\partial \Sigma_t = \Gamma_0 \cup \Gamma_t$ and $\chi(\Sigma_t) = 0$ is the Euler characteristic of an annulus. Differentiation of \eqref{eq:gauss_bonnet} and subsequent substitution yields
    \begin{align*}
        \frac{\diff}{\diff t} \int_{\Gamma_t} \psi_1 \diff s = - \int_{\Gamma_t} \kappa K \diff s,
    \end{align*}
    where $\psi_1$ is the geodesic curvature of $\Gamma_t$ on $\Sigma_t$ and the integrand on the right hand side is $0$ by the assumption that $K = 0$. 
\end{proof}
\end{proposition}

\noindent As in the previous subsection, we construct analytical examples using configurations with helical and cylindrical symmetries.

\begin{example}[Helical symmetry]\label{ex:cgc_helical_symmetry}
    For a constant $\theta_0$ consider evolving helix
    \begin{align*}
        \gamma_0(u) &:= \begin{bmatrix}  
        \varrho_0 \cos u \\
        \varrho_0 \sin u \\
        wu
        \end{bmatrix}, &
        \gamma(t, u) &:= \begin{bmatrix}  
        \varrho(t) \cos(u + \upsilon(t)) \\
        \varrho(t) \sin(u + \upsilon(t)) \\
        wu + \omega(t)
        \end{bmatrix},
    \end{align*}
    where $\rho_0$ and $w$ are positive constants and $\rho, \upsilon, \omega$ are functions of time $t$. Since $\kappa = \varrho g^{-2}$ and $\tau = w g^{-2}$, the problem \eqref{eq:cgc-motion-law} reduces to the following system
    \begin{align}\label{eq:cgc_helical_example}
        \frac{\diff}{\diff t} \begin{bmatrix} \theta \\ \varrho \\ \omega \\ \upsilon \end{bmatrix} = 
        \frac{1}{g^3}
        \begin{bmatrix}
            \frac{K g^4 + w^2 \cos^2\theta}{g \sin \theta} \\
            - g \rho \cos\theta \\
            \rho^2 \sin\theta \\
            - w \sin\theta
        \end{bmatrix}, &&
        \left.
        \begin{bmatrix} \theta \\ \varrho \\ \omega \\ \upsilon \end{bmatrix}
        \right\vert_{t = 0}=
        \begin{bmatrix} \theta_0 \\ \varrho_0 \\ 0 \\ 0 \end{bmatrix},
    \end{align}
    where $g^2 = \varrho^2 + w^2$. This solution leads to a family of helical trajectory surfaces of constant Gaussian curvature.
\end{example}

\begin{example}[Cylindrical symmetry]\label{ex:cylindrical_symmetry_gauss} 
    Setting $w = 0$ reduces \eqref{eq:cgc_helical_example} to
    \begin{align*}
        \frac{\diff}{\diff t} \begin{bmatrix} \theta \\ \varrho \\ \omega \end{bmatrix} = 
        \frac{1}{\varrho^2}
        \begin{bmatrix}
            (\sin\theta )^{-1}K g^3 \\
            - \varrho \cos\theta \\
            \varrho \sin\theta
        \end{bmatrix}, &&
        \left.
        \begin{bmatrix} \theta \\ \varrho \\ \omega \end{bmatrix}
        \right\vert_{t = 0}=
        \begin{bmatrix} \theta_0 \\ \varrho_0 \\ 0 \end{bmatrix}.
    \end{align*}
\end{example}

%%%%%%%%%%%%%%%%%%%%%%%%%%%%%%%%%%%%%%%%%%%%%%%%%%%%%%%%%%%%%%%%%%%%%%%
\section{Conclusion}

The framed curvature flow offers a promising extension of classical curvature-driven geometric flows of space curves, providing a rich configuration space of motion laws with potential applications in the analysis of surfaces with prescribed curvature and beyond.

In this work, we introduced the framed curvature flow as a generalization of the curve shortening flow and vortex filament equation, incorporating a time-dependent moving frame to govern the direction of curvature-driven motion. We established local existence and uniqueness for a simplified version of this flow, derived global estimates for key geometric and topological quantities, and examined the trajectory surfaces produced by variations of the flow, including those leading to surfaces of constant mean or Gaussian curvature. Additionally, we classified singularities that may arise during the flow.

More research is needed to explore local existence across a broader range of $\theta$-velocity settings. Future work could also investigate motion laws that generate other surface types, such as those with a constant principal curvature ratio \cite{new_pottman}, or surfaces that minimize energies like Willmore energy \cite{soliman_willmore} or various repulsive energies \cite{repulsive_surfaces,repulsive_curves}.

Another direction for future research could involve extending the concept of framed curvature flow to higher-dimensional spaces with more than one codimension, and exploring possible connections between the generated trajectory varieties and Open book decomposition.

Further insights may be gained through rigorous numerical analysis and experiments with different $\theta$-velocity settings. Finally, the examples of singularities involving curvature blow-up, discussed in Subsection~\ref{subsec:analysis_of_singularities}, should be expanded and studied in greater detail.

%%%%%%%%%%%%%%%%%%%%%%%%%%%%%%%%%%%%%%%%%%%%%%%%%%%%%%%%%%%%%%%%%%%%%%%

\end{document}